\documentclass{amsart}

\usepackage{amscd,amssymb,amsmath,graphicx,verbatim}
\usepackage[dvips]{hyperref}
\usepackage{textcomp}
\usepackage[OT1,T1]{fontenc}

\theoremstyle{definition}

\theoremstyle{remark}

\numberwithin{equation}{section}

\begin{document}

\title[An additional reformulation of the RH]{Other representations of the Riemann Zeta function and an
additional reformulation of the Riemann Hypothesis}
\author{Stefano Beltraminelli}
\address{S. Beltraminelli, CERFIM, Research Center for Mathematics
and Physics, PO Box 1132, 6600 Locarno, Switzerland}
\email{stefano.beltraminelli@ti.ch}
\author{Danilo Merlini}
\address{D. Merlini, CERFIM, Research Center for Mathematics and
Physics, PO Box 1132, 6600 Locarno, Switzerland}
\email{merlini@cerfim.ch}
\label{I1}
\date{15 April 2007}
\subjclass{11M26}
\keywords{Riemann Zeta function, Riemann Hypothesis, Criteria of
Riesz, Hardy-Littlewood and Baez-Duarte, Pochhammer's polynomials}
\begin{abstract}
New expansions for some functions related to the Zeta function in
terms of the Pochhammer's polynomials are given (coefficients $b_{k}$,
$d_{k}$, ${\hat{d}}_{k}$ and ${\hat{\hat{d}}}_{k}$). In some formal
limit our expansion $b_{k}$ obtained via the alternating series
gives the regularized expansion of Maslanka for the Zeta function.
The real and the imaginary part of the function on the critical
line is obtained with a good accuracy up to $\mathfrak{I}( s) =t<35$.

Then, we give the expansion (coefficient $\hat{d_{k}}$) for the
derivative of $\ln ( (s-1)\zeta ( s) ) $. The critical function
of the derivative, whose bounded values for\ \ $\mathfrak{R}( s)
>\frac{1}{2}$ at large values of {\itshape k} should ensure the
truth of the Riemann Hypothesis (RH), is obtained either by means
of the primes or by means of the zeros (trivial and non-trivial)
of the Zeta function. In a numerical experiment performed up to
high values of {\itshape k} i.e. up to $k={10}^{13}$ we obtain a
very good agreement between the two functions, with the emergence
of twelve oscillations with stable amplitude.

For a special case of values of the two parameters entering in the
general Pochhammer's expansion it is argued that the bound on the
critical function should be given by the Euler constant gamma.
\end{abstract}
\maketitle

\section{Introduction}

Lately there has been new interest in the study of the expansion
of the Zeta function via the Pochhammer's polynomials. This is related
to the original idea of Riesz \cite{11} and of Hardy-Littlewood
\cite{12} at the beginning of the last century. In a pioneering
work \cite{10} Maslanka obtained a regularized expansion for the
Zeta function (with coefficients $A_{k}$) and Baez-Duarte for the
expansion of the reciprocal of the Zeta function (with coefficients
$c_{k}$) for the Riesz case \cite{2,3}. Other cases of interest
have also recently been studied \cite{6,7,8,9,1}. As pointed out
in \cite{3}, the discrete version by means of the Pochhammer's polynomials
$P_{k}( s) $, where $s=\sigma +i t$ is the complex variable and
{\itshape k} is an integer, has advantages especially in the context
of numerical experiments in connection with some ``kind of verification''
in the direction to believe that the RH may be true. 

In this work we first derive a new expansion for the Zeta function
in terms of the Pochhammer's polynomials via the alternating series
(with new coefficients $b_{k}$). In some formal limit, a connection
with the expansion of Maslanka is also obtained in Section 2. Our
expansion is then studied numerically on the critical line where
a good agreement with the real function is obtained up to $\mathfrak{I}(
s) =t<35$, with the emergence of the first few low zeros. After
this value of {\itshape t}, a divergence possibly of numerical nature
set on.

In Section 3 we then obtain the expansion for the function $\ln
( (1-2^{1-s})\zeta ( s) ) $ (with new coefficients $d_{k}$) as well
for the derivative of $\ln ( (s-1)\zeta ( s) ) $ (with new coefficients
$\hat{d_{k}}$) in terms of the two parameters $\alpha$ and $\beta$,
already introduced in our previous works \cite{13,5,14} Then the
critical function for the derivative (whose boundedness at large
{\itshape k} would ensure the truth of the RH) is then obtained
either with the primes or with the trivial and non-trivial zeros
of the Zeta function.

In a numerical experiment for the special case $\alpha =\frac{9}{2}$
and $\beta =4$ up to high values of {\itshape k}, i.e. $k={10}^{13}$,
the results for the two functions are in very good agreement, both
with the emergence of the same twelve oscillations of stable amplitude
of about 0.01 (Section 4).

Finally, in the limit of large $\beta$ and $\alpha =1$, it is argued
that an upper bound to the critical function should be given by
the Euler constant gamma (Section 5).

\section{Zeta function representation via the alternating series}

In this section we derive a formula for $(1-2^{1-s})\zeta ( s) $
similar to the one of Maslanka \cite{10} for $(s-1)\zeta ( s) $
and of Baez-Duarte \cite{2,3} for ${[\zeta ( s) ]}^{-1}$.

Here the starting series is convergent for $\mathfrak{R}( s) =\sigma
>0$ and the formula is obtained still in terms of the so called
Pochhammer's polynomials of degree {\itshape k}, in the complex
variable $s=\sigma +i t$.
\begin{equation}
P_{k}( s) =\prod \limits_{r=1}^{k}\left( 1-\frac{s}{r}\right) \ \ \ \ \ \ \ \ \ \ \forall
k\in \mathbb{N}^{*}\ \ \ \ \mathrm{and}\ \ \ P_{0}( s) =1 
\end{equation}

We will also use a family of functions with two parameters ($\alpha$
and $\beta$) as considered already in our recent works \cite{5,13,14}.
Since the alternating series is given by:
\begin{equation}
\left( 1-2^{1-s}\right) \zeta ( s) =\sum \limits_{n=1}^{\infty }\frac{{\left(
-1\right) }^{n-1}}{n^{s}}\ \ \ \ \ \ \ \ \forall \mathfrak{R}( s)
=\sigma >0
\end{equation}

\noindent we have using the trick as in \cite{2} that:
\[
\begin{array}{rl}
 \left( 1-2^{1-s}\right) \zeta ( s)  & =\sum \limits_{n=1}^{\infty
}\frac{{\left( -1\right) }^{n-1}}{n^{\alpha }}{\left( 1-\left( 1-\frac{1}{n^{\beta
}}\right) \right) }^{\frac{s-\alpha }{\beta }} \\
  & =\sum \limits_{n=1}^{\infty }\frac{{\left( -1\right) }^{n-1}}{n^{\alpha
}}\sum \limits_{k=0}^{\infty }{\left( -1\right) }^{k}{\left( 1-\frac{1}{n^{\beta
}}\right) }^{k}\binom{\frac{s-\alpha }{\beta }}{k}
\end{array}
\]

Since
\[
\begin{array}{rl}
 {\left( -1\right) }^{k}\binom{\frac{s-\alpha }{\beta }}{k} & =\frac{{\left(
-1\right) }^{k}}{k!}\left( \frac{s-\alpha }{\beta }+1-1\right) \cdots
( \frac{s-\alpha }{\beta }+1-k)  \\
  & =\prod \limits_{r=1}^{k}\left( 1-\frac{\frac{s-\alpha }{\beta
}+1}{r}\right) =P_{k}( \frac{s-\alpha }{\beta }+1) 
\end{array}
\]

\noindent we obtain:
\begin{equation}
\begin{array}{rl}
 \left( 1-2^{1-s}\right) \zeta ( s)  & =\sum \limits_{k=0}^{\infty
}P_{k}( \frac{s-\alpha }{\beta }+1) \sum \limits_{n=1}^{\infty }\frac{{\left(
-1\right) }^{n-1}}{n^{\alpha }}{\left( 1-\frac{1}{n^{\beta }}\right)
}^{k} \\
  & =\sum \limits_{k=0}^{\infty }P_{k}( \frac{s-\alpha }{\beta }+1)
\sum \limits_{j=0}^{k}{\left( -1\right) }^{j}\binom{k}{j}\sum \limits_{n=1}^{\infty
}\frac{{\left( -1\right) }^{n-1}}{n^{\alpha +\beta j}}
\end{array}
\end{equation}

Since from (2.2)
\[
\left( 1-2^{1-\left( \alpha +\beta j\right) }\right) \zeta ( \alpha
+\beta j) =\sum \limits_{n=1}^{\infty }\frac{{\left( -1\right) }^{n-1}}{n^{\alpha
+\beta j}}
\]

\noindent substitution in (2.3) gives:
\begin{equation}
\left( 1-2^{1-s}\right) \zeta ( s) =\sum \limits_{k=0}^{\infty }P_{k}(
\frac{s-\alpha }{\beta }+1) \sum \limits_{j=0}^{k}{\left( -1\right)
}^{j}\binom{k}{j}\left( 1-2^{1-\left( \alpha +\beta j\right) }\right)
\zeta ( \alpha +\beta j\mathrm{)} 
\end{equation}

With the definition
\begin{equation}
b_{k}:=\sum \limits_{j=0}^{k}{\left( -1\right) }^{j}\binom{k}{j}\left(
1-2^{1-\left( \alpha +\beta j\right) }\right) \zeta ( \alpha +\beta
j\mathrm{)} \mathrm{\ \ }
\end{equation}

\noindent (2.4) becomes:
\begin{equation}
\left( 1-2^{1-s}\right) \zeta ( s) =\sum \limits_{k=0}^{\infty }b_{k}P_{k}(
\frac{s-\alpha }{\beta }+1) 
\end{equation}

\noindent where $P_{0}( \frac{s-\alpha }{\beta }+1) =1$ and $b_{0}=(1-2^{1-\alpha
})\zeta ( \alpha ) $.

The series above, is expected to represent $(1-2^{1-s})\zeta ( s)
$ for {\itshape s} in some compact subset of the plane as for the
Maslanka case \cite{10}. In that case, the central point has been
investigated and elucidated by Baez-Duarte \cite{15}. Here many
choices of $\alpha$ and $\beta$ are possible. For $\alpha =\beta
=2$ we have the Riesz case \cite{11} and it is the analogon to the
regularized version of Maslanka but the representation of the Zeta
function is not the same. For $\alpha =1+\delta \ \ (\delta \downarrow
0)$ and $\beta =2$ we obtain the Hardy-Littlewood case \cite{12}
which was also discussed numerically in a different way using other
polynomials \cite{17}.

In fact, from Lemma 2.3 of Baez-Duarte \cite{3} which states that
at large {\itshape k}:
\begin{equation}
\left| P_{k}( s) \right| \leq C k^{-\mathfrak{R}( s) }
\end{equation}

\noindent  where C is a constant depending on $ |$s$ |$, 

we obtain here that:
\[
\left| P_{k}( \frac{s-\alpha }{\beta }+1) \right| \leq C k^{-\left(
\frac{\mathfrak{R}( s) -\alpha }{\beta }+1\right) }
\]

We thus suspect and expect that the above series represents $(1-2^{1-s})\zeta
( s) $ for all $\mathfrak{R}( s) >\frac{1}{2}+\delta , \delta >0$
if we assume $|b_{k}|\leq D k^{-\gamma }$ with $\gamma \geq \frac{\alpha
-1/2-\delta }{\beta }$ at large values of {\itshape k }and for some
constant {\itshape D. }In fact with this assumption we have that:
\[
\begin{array}{rl}
 \left| \left( 1-2^{1-s}\right) \zeta ( s) \right|  & \leq \sum
\limits_{k=0}^{\infty }\left| b_{k}P_{k}( \frac{s-\alpha }{\beta
}+1) \right| \leq \mathrm{const}.\sum \limits_{k=0}^{\infty }k^{-\frac{\alpha
-1/2-\delta }{\beta }}k^{-\left( \frac{\mathfrak{R}( s) -\alpha
}{\beta }+1\right) } \\
  & \leq \mathrm{const}.\sum \limits_{k=0}^{\infty }k^{-\left( 1+\frac{\mathfrak{R}(
s) -1/2-\delta }{\beta }\right) }<\infty 
\end{array}
\]

\noindent if $\mathfrak{R}( s) >\frac{1}{2}+\delta $. 

For $\alpha =\beta =2$ (Riesz) we should have $|b_{k}|\leq D k^{-\frac{3}{4}+\epsilon
}$. For the case $\alpha =1$ and $\beta =2$ (Hardy-Littlewood) we
should have $|b_{k}|\leq D k^{-\frac{1}{4}+\epsilon }$. Another
case of interest is the one where $\alpha =\frac{3}{2}$ and $\beta
=1$. In this case one should have $|b_{k}|\leq D k^{-1+\epsilon
}$.\ \ \ \ 

Of interest also, is the limiting case of large values of $\beta$,
where barely $b_{k}$ should behave as $|b_{k}|\leq D$. 

For a strong argument (a Theorem) in favour of the validity of the
Maslanka representation of $(s-1)\zeta ( s) $ in some regions of
the complex plane (compact subsets), the reader should consult the
works of Baez-Duarte \cite{15} already mentionned and it is expected
that using the same methods, the proof of (2.6) may be obtained
for all $\mathfrak{R}( s) >\frac{1}{2}$. Here, for our series we
limit ourselves to a numerical analysis just illustrating the kind
of accuracy of some representations.

{\itshape Remark}. Let us consider the Riesz case $\alpha =\beta
=2$. We can write:
\[
\left( 1-e^{\left( 1-s\right) \ln  2}\right) \zeta ( s) =\sum \limits_{k=0}^{\infty
}P_{k}( \frac{s}{2}) \sum \limits_{j=0}^{k}{\left( -1\right) }^{j}\binom{k}{j}\left(
1-e^{-\left( 1+2j\right) \ln  \mathrm{2}}\right) \zeta ( 2+2j\mathrm{)}
\]

\noindent and using the Taylor's expansion of $ e ^{x}$, we obtain:
\begin{equation}
\left( s-1\right) \zeta ( s) =\sum \limits_{k=0}^{\infty }A_{k}P_{k}(
\frac{s}{2}) 
\end{equation}

\noindent where
\begin{equation}
A_{k}=\sum \limits_{j=0}^{k}{\left( -1\right) }^{j}\binom{k}{j}\left(
2j+1\right) \zeta ( 2j+2) 
\end{equation}

\noindent i.e. the representation obtained originally by a different
method by Maslanka in a pioneering work \cite{10}. We remark that
(2.8) and (2.9) should not be considered as an approximation of
our formulas (2.5) and (2.6) and vice versa. (2.5), (2.6) and (2.8),
(2.9) are simply two different representations of functions related
to the Riemann Zeta function, the first one given by $(s-1)\zeta
( s) $, the second one by $(1-2^{1-s})\zeta ( s) $.

As an example, for $s=\sigma $ with $\sigma$ in $[0,1]$, both representations
give a good description of the real function $\zeta ( \sigma ) $
as may easily be computationally checked.

We now proceed to obtain a representation of $\zeta ( s) $ possibly
correct on the critical line $s=\frac{1}{2}+i t$, with the help
of (2.5) and (2.6), in which we are free to set $\alpha =\frac{1}{2}$
and $\beta =i$. Then:
\begin{equation}
\left( 1-2^{\frac{1}{2}-i t}\right) \zeta ( \frac{1}{2}+i t) =\sum
\limits_{k=0}^{\infty }b_{k}P_{k}( t+1) 
\end{equation}

\noindent where now
\begin{equation}
b_{k}=\sum \limits_{j=0}^{k}{\left( -1\right) }^{j}\binom{k}{j}\left(
1-2^{\frac{1}{2}-ij}\right) \zeta ( \frac{1}{2}+ij\mathrm{)} \mathrm{\ \ }
\end{equation}

We now check the series in (2.10) restricting {\itshape k} up to
20 for $t\leq 18$ and up to 50 for $t>18$. We compare the result
with the exact functions $\mathfrak{R}( (1-2^{\overline{s}})\zeta
( s) ) $ and $\mathfrak{I}( (1-2^{\overline{s}})\zeta ( s) ) $,
for $s=\frac{1}{2}+i t$ with {\itshape t} up to 40. The plots are
given below. We obtain a good approximation with the emergence of
the first five non-trivial zeros located at $t_{1}=14.13472\ldots
$, $t_{2}=21.02204\ldots $, $t_{3}=25.01085\ldots $, $t_{4}=21.02204\ldots
$, $t_{5}=32.93505\ldots $. The numerical results are satisfactory
until $t\cong 35$.\ \ 
\begin{figure}[h]
\begin{center}
\includegraphics{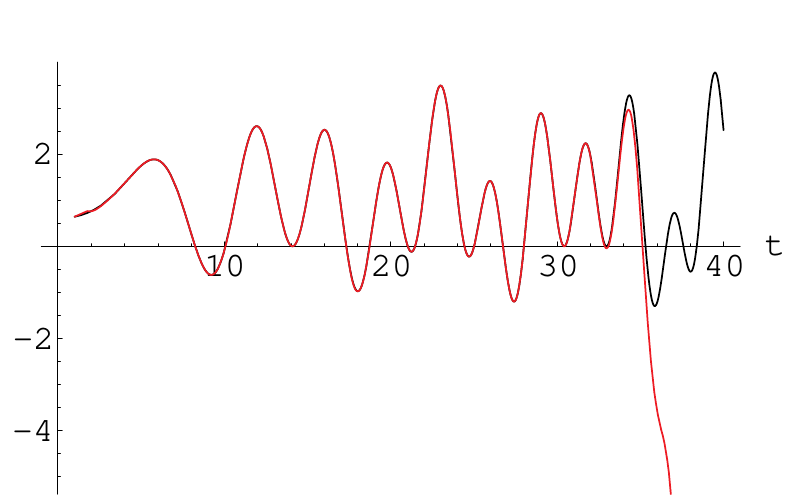}

\end{center}
\caption{The plot of the real part of $\sum \limits_{k=0}^{20 (50)}b_{k}P_{k}(
t+1) $ [red] vs. $\mathfrak{R}( (1-2^{\overline{s}})\zeta ( s) )
$ [black]}

\end{figure}
\begin{figure}[h]
\begin{center}
\includegraphics{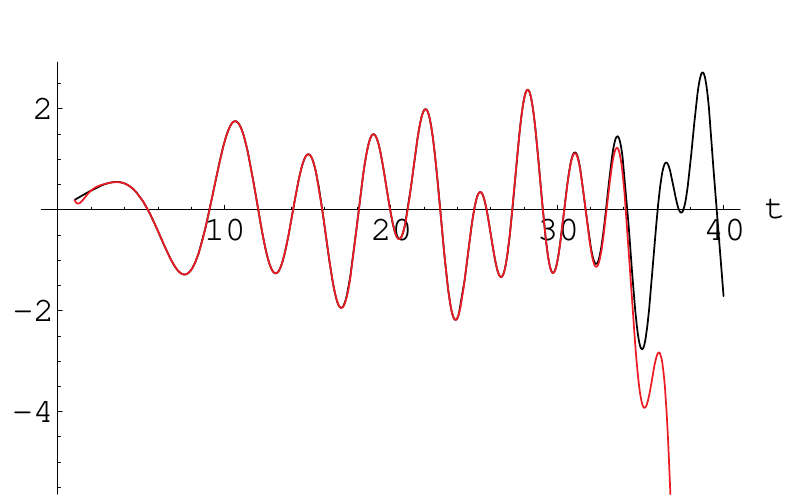}

\end{center}
\caption{The plot of the imaginary part of $\sum \limits_{k=0}^{20
(50)}b_{k}P_{k}( t+1) $ [red] vs. $\mathfrak{I}( (1-2^{\overline{s}})\zeta
( s) ) $ [black]}

\end{figure}

This concludes the first part of our work. Below, in the second
part we develop two new representations of the functions $\ln (
(1-2^{1-s})\zeta ( s) ) $ and $\frac{d}{ds}\ln ( (s-1)\zeta ( s)
) $ which may possibly constitute a satisfactory approximation to
the exact functions.

\section{A representation for the logarithm of the Zeta Function
and an additional criterion for the truth of the RH}

We will start as before but instead of writing $\zeta ( s) $ as
a sum, i.e. $\zeta ( s) =\sum \limits_{n=1}^{\infty }\frac{1}{n^{s}}$,
we will use the Euler product formula to derive a new representation
for $\ln ( (1-2^{1-s})\zeta ( s) ) $, which of course should be
carefully investigated by means of some numerical experiments. Thus:
\begin{equation}
\ln [ \left( 1-2^{1-s}\right) \zeta ( s) ] =\ln [ \left( 1-2^{1-s}\right)
\prod \limits_{p\,\mathrm{prime}}\frac{1}{1-p^{-s}}] \ \ \ \ \ \ \ \forall
\mathfrak{R}( s) >1
\end{equation}

For any prime {\itshape p}, we have:
\[
\ln ( 1-p^{-s}) =-\sum \limits_{n=1}^{\infty }\frac{p^{-\mathit{ns}}}{n}
\]

\noindent so that introducing the parameters $\alpha$ and $\beta$
as before we have that:
\[
\begin{array}{rl}
 \sum \limits_{n=1}^{\infty }\frac{p^{-\alpha n}}{n}{\left( 1-\left(
1-p^{-\beta n}\right) \right) }^{\frac{s-\alpha }{\beta }} & =\sum
\limits_{n=1}^{\infty }\frac{p^{-\alpha n}}{n}\sum \limits_{k=0}^{\infty
}{\left( -1\right) }^{k}{\left( 1-p^{-\beta n}\right) }^{k}\binom{\frac{s-\alpha
}{\beta }}{k} \\
  & =\sum \limits_{k=0}^{\infty }P_{k}( \frac{s-\alpha }{\beta }+1)
\sum \limits_{n=1}^{\infty }\frac{1}{n}\sum \limits_{j=0}^{k}{\left(
-1\right) }^{j}\binom{k}{j}p^{-\left( \alpha +\beta j\right) n}
 \\
  & =\sum \limits_{k=0}^{\infty }P_{k}( \frac{s-\alpha }{\beta }+1)
\sum \limits_{j=0}^{k}{\left( -1\right) }^{j}\binom{k}{j}\ln ( 1-p^{-\left(
\alpha +\beta j\right) }) 
\end{array}
\]

\noindent the same treatment for the function $\ln ( 1-2^{1-s})
$, gives:
\[
\ln ( 1-2^{1-s}) =\sum \limits_{k=0}^{\infty }P_{k}( \frac{s-\alpha
}{\beta }+1) \sum \limits_{j=0}^{k}{\left( -1\right) }^{j}\binom{k}{j}\ln
( 1-2^{1-\left( \alpha +\beta j\right) }) 
\]

\noindent where $P_{k}$ are still the Pochhammer's polynomials.

Finally, the representation of $\ln ( (1-2^{1-s})\zeta ( s) ) $,
we propose is given by:
\begin{equation}
\ln [ \left( 1-2^{1-s}\right) \zeta ( s) ] =\sum \limits_{k=0}^{\infty
}d_{k}P_{k}( \frac{s-\alpha }{\beta }+1) 
\end{equation}

\noindent where now:
\begin{equation}
d_{k}:=\sum \limits_{j=0}^{k}{\left( -1\right) }^{j}\binom{k}{j}\ln
[ \left( 1-2^{1-\left( \alpha +\beta j\right) }\right) \zeta ( \alpha
+\beta j) ] 
\end{equation}

{\itshape Remark}. Another formal derivation of the above equations
is the following:
\[
\ln [ \left( 1-2^{1-s}\right) \zeta ( s) ] =\ln \left[ \sum \limits_{n=1}^{\infty
}\frac{{\left( -1\right) }^{n-1}}{n^{s}}\right] 
\]

Supposing now that the right hand side may be given as an unknown
series $\sum \limits_{r=1}^{\infty }\frac{a_{r}}{r^{s}}$ we then
have:
\[
\begin{array}{rl}
 \sum \limits_{r=1}^{\infty }\frac{a_{r}}{r^{\alpha }}{\left( 1-\left(
1-\frac{1}{r^{\beta }}\right) \right) }^{\frac{s-\alpha }{\beta
}} & =\sum \limits_{k=0}^{\infty }P_{k}( \frac{s-\alpha }{\beta
}+1) \sum \limits_{r=1}^{\infty }\frac{a_{r}}{r^{\alpha }}{\left(
1-\frac{1}{r^{\beta }}\right) }^{k} \\
  & =\sum \limits_{k=0}^{\infty }P_{k}( \frac{s-\alpha }{\beta }+1)
\sum \limits_{j=0}^{k}{\left( -1\right) }^{j}\binom{k}{j}\sum \limits_{r=1}^{\infty
}\frac{a_{r}}{r^{\alpha +\beta j}} \\
  & =\sum \limits_{k=0}^{\infty }P_{k}( \frac{s-\alpha }{\beta }+1)
\sum \limits_{j=0}^{k}{\left( -1\right) }^{j}\binom{k}{j}\ln \left(
\sum \limits_{n=1}^{\infty }\frac{{\left( -1\right) }^{n-1}}{n^{\alpha
+\beta j}}\right) 
\end{array}
\]

\noindent which coincide with (3.2) and (3.3), obtained with the
Euler product formula for $\mathfrak{R}( s) >1$. (3.2) with (3.3),
is the new formula possibly representing the logarithm of the Zeta
function in terms of the two parameters Pochhammer's polynomials.
To the best of our knowledge the above representation is new and
it is our aim to carry out some numerical experiments in the sequel
in order to support its validity also in some compact subset of
the critical strip.

We now investigate the representation of the derivative of the function
$\mathit{\ln }( (s-1)\zeta ( s) ) $: 
\begin{equation}
\frac{d}{ds}\ln ( \left( s-1\right) \zeta ( s) ) =\frac{1}{s-1}+\frac{\zeta
^{\prime }( s) }{\zeta ( s) }
\end{equation}

Then with $\zeta ( s) =\prod \limits_{p\,\mathrm{prime}}\frac{1}{1-p^{-s}}$
we obtain:
\[
\begin{array}{rl}
 \frac{\zeta ^{\prime }( s) }{\zeta ( s) } & =-\sum \limits_{p }\frac{d}{ds}\ln
( 1-p^{-s}) =-\sum \limits_{p}\frac{1}{1-p^{-s}}\frac{d}{ds}\left(
1- e ^{-s \ln  p}\right)  \\
  & =-\sum \limits_{p}\frac{p^{-s}}{1-p^{-s}}\ln  p=-\sum \limits_{p}\ln
p\sum \limits_{q=1}^{\infty }\frac{1}{p^{sq}}
\end{array}
\]

Introducing as above the Pochhammer's polynomials we obtain further:
\[
\begin{array}{rl}
 \frac{\zeta ^{\prime }( s) }{\zeta ( s) } & =-\sum \limits_{p}\ln
p\sum \limits_{q=1}^{\infty }\frac{1}{p^{q\alpha }}{\left( 1-\left(
1-\frac{1}{p^{q\beta }}\right) \right) }^{\frac{s-\alpha }{\beta
}} \\
  & =-\sum \limits_{p}\ln  p\sum \limits_{k=0}^{\infty }P_{k}( \frac{s-\alpha
}{\beta }+1) \sum \limits_{j=0}^{k}{\left( -1\right) }^{j}\binom{k}{j}\sum
\limits_{q=1}^{\infty }\frac{1}{p^{q( \alpha +\beta j) }} \\
  & =\sum \limits_{k=0}^{\infty }P_{k}( \frac{s-\alpha }{\beta }+1)
\sum \limits_{j=0}^{k}{\left( -1\right) }^{j}\binom{k}{j}\sum \limits_{q=1}^{\infty
}\left( -\sum \limits_{p}\frac{1}{p^{q( \alpha +\beta j) }}\ln 
p\right)  \\
  & =\sum \limits_{k=0}^{\infty }P_{k}( \frac{s-\alpha }{\beta }+1)
\sum \limits_{j=0}^{k}{\left( -1\right) }^{j}\binom{k}{j}\frac{\partial
}{\partial \alpha }\left( \sum \limits_{q=1}^{\infty }\frac{1}{q}\sum
\limits_{p}\frac{1}{p^{q( \alpha +\beta j) }}\right)  \\
  & =\sum \limits_{k=0}^{\infty }P_{k}( \frac{s-\alpha }{\beta }+1)
\sum \limits_{j=0}^{k}{\left( -1\right) }^{j}\binom{k}{j}\frac{\partial
}{\partial \alpha }\left( -\sum \limits_{p}\ln ( 1-\frac{1}{p^{\alpha
+\beta j}}) \right)  \\
  & =\sum \limits_{k=0}^{\infty }P_{k}( \frac{s-\alpha }{\beta }+1)
\sum \limits_{j=0}^{k}{\left( -1\right) }^{j}\binom{k}{j}\frac{\partial
}{\partial \alpha }\ln ( \prod \limits_{p}\frac{1}{1-p^{-\left(
\alpha +\beta j\right) }})  \\
  & =\sum \limits_{k=0}^{\infty }P_{k}( \frac{s-\alpha }{\beta }+1)
\sum \limits_{j=0}^{k}{\left( -1\right) }^{j}\binom{k}{j}\frac{\partial
}{\partial \alpha }\ln  \zeta ( \alpha +\beta j) 
\end{array}
\]

For $\frac{1}{s-1}$, using $\frac{1}{s-1}$=$\int _{0}^{\infty }
e ^{-\lambda ( s-1) }d\lambda $ we have similarly: 
\[
\begin{array}{rl}
 \frac{1}{s-1} & =\int _{0}^{\infty } e ^{\lambda }\frac{1}{ e ^{\lambda
s}}d\lambda =\int _{0}^{\infty }\frac{ e ^{\lambda }}{ e ^{\lambda
\alpha }}{\left( 1-\left( 1-\frac{1}{ e ^{\lambda \beta }}\right)
\right) }^{\frac{s-\alpha }{\beta }}d\lambda  \\
  & =\int _{0}^{\infty } e ^{\lambda }\sum \limits_{k=0}^{\infty
}P_{k}( \frac{s-\alpha }{\beta }+1) \sum \limits_{j=0}^{k}{\left(
-1\right) }^{j}\binom{k}{j}\frac{1}{ e ^{\lambda ( \alpha +\beta
j) }}d\lambda  \\
  & =\sum \limits_{k=0}^{\infty }P_{k}( \frac{s-\alpha }{\beta }+1)
\sum \limits_{j=0}^{k}{\left( -1\right) }^{j}\binom{k}{j}\int _{0}^{\infty
} e ^{-\lambda ( \alpha +\beta j-1) }d\lambda  \\
  & =\sum \limits_{k=0}^{\infty }P_{k}( \frac{s-\alpha }{\beta }+1)
\sum \limits_{j=0}^{k}{\left( -1\right) }^{j}\binom{k}{j}\frac{1}{\alpha
+\beta j-1} \\
  & =\sum \limits_{k=0}^{\infty }P_{k}( \frac{s-\alpha }{\beta }+1)
\sum \limits_{j=0}^{k}{\left( -1\right) }^{j}\binom{k}{j}\frac{\partial
}{\partial \alpha }\ln ( \alpha +\beta j-1) 
\end{array}
\]

Thus, along these lines we obtain:
\begin{equation}
\frac{d}{ds}\ln ( \left( s-1\right) \zeta ( s) ) =\sum \limits_{k=0}^{\infty
}{\hat{d}}_{k}P_{k}( \frac{s-\alpha }{\beta }+1) 
\end{equation}

\noindent where:
\begin{equation}
{\hat{d}}_{k}=\sum \limits_{j=0}^{k}{\left( -1\right) }^{j}\binom{k}{j}\frac{\partial
}{\partial \alpha }\ln [ \left( \alpha +\beta j-1\right) \zeta (
\alpha +\beta j) ] 
\end{equation}

From the formula (7) in \cite{18}, where $\rho$ represents a non-trivial
zero of the Zeta function, i.e.:
\[
\begin{array}{rl}
 \frac{1}{s-1}+\frac{\zeta ^{\prime }( s) }{\zeta ( s) } & =\frac{1}{s-1}-\frac{s}{s-1}+\sum
\limits_{\rho }\frac{1}{\rho }+\sum \limits_{\rho }\frac{1}{s-\rho
}-\sum \limits_{n=1}^{\infty }\frac{1}{2n}+\sum \limits_{n=1}^{\infty
}\frac{1}{s+2n}+\frac{\zeta ^{\prime }( 0) }{\zeta ( 0) } \\
  & =\frac{\zeta ^{\prime }( 0) }{\zeta ( 0) }-1+\sum \limits_{\rho
}\frac{1}{\rho }-\sum \limits_{n=1}^{\infty }\frac{1}{2n}+\sum \limits_{\rho
}\frac{1}{s-\rho }+\sum \limits_{n=1}^{\infty }\frac{1}{s+2n}
\end{array}
\]

Setting $C=\frac{\zeta ^{\prime }( 0) }{\zeta ( 0) }-1$, this equation
applied to $s=\alpha +\beta j$ in (3.6) gives:
\[
\begin{array}{rl}
 {\hat{d}}_{k} & =\sum \limits_{j=0}^{k}{\left( -1\right) }^{j}\binom{k}{j}\left(
C+\int _{0}^{\infty }\left( \sum \limits_{\rho } e ^{-\lambda (
\alpha +\beta j-\rho ) }+ e ^{-\lambda \rho }+\sum \limits_{n=1}^{\infty
} e ^{-\lambda ( \alpha +\beta j+2n) }- e ^{-\lambda 2n}\right)
d\lambda \right)  \\
  & \begin{array}{l}
 =\int _{0}^{\infty }\sum \limits_{\rho }\left( { e ^{-\lambda (
\alpha -\rho ) }( 1-\frac{1}{ e ^{\lambda \beta }}) }^{k}+ e ^{-\lambda
}{\left( 1-\frac{1}{ e ^{\lambda \beta }}\right) }^{k}\delta _{k,0}\right)
d\lambda + \\
 \int _{0}^{\infty }\left( \sum \limits_{n=1}^{\infty }{ e ^{-\lambda
( \alpha +2n) }( 1-\frac{1}{ e ^{\lambda \beta }}) }^{k}- e ^{-\lambda
2n}{\left( 1-\frac{1}{ e ^{\lambda \beta }}\right) }^{k}\delta _{k,0}\right)
d\lambda 
\end{array}
\end{array}
\]

We consider only $k>0$. Now we make the variable change $ e ^{-\lambda
\beta }=x$ and finally we obtain: 
\[
\begin{array}{rl}
 {\hat{d}}_{k} & =\frac{1}{\beta }\int _{0}^{1}{\left( 1-x\right)
}^{k+1-1}\sum \limits_{\rho }x^{\frac{\alpha -\rho }{\beta }-1}dx+\frac{1}{\beta
}\int _{0}^{1}{\left( 1-x\right) }^{k+1-1}\sum \limits_{n=1}^{\infty
}x^{\frac{\alpha +2n}{\beta }-1}dx \\
  & =\frac{1}{\beta }\sum \limits_{\rho }B( \frac{\alpha -\rho }{\beta
},k+1) +\frac{1}{\beta }\sum \limits_{n=1}^{\infty }B( \frac{\alpha
+2n}{\beta },k+1) 
\end{array}
\]

\noindent where $B( x,y) =\frac{\Gamma ( x) \Gamma ( y) }{\Gamma
( x+y) }$ is the Beta function. Thus for large {\itshape k} we can
write:
\begin{equation}
{\hat{d}}_{k}=\frac{1}{\beta }\sum \limits_{\rho }\Gamma ( \frac{\alpha
-\rho }{\beta }) k^{-\frac{\alpha -\rho }{\beta }}+\frac{1}{\beta
}\sum \limits_{n=1}^{\infty }\Gamma ( \frac{\alpha +2n}{\beta })
k^{-\frac{\alpha +2n}{\beta }}
\end{equation}

For the critical function \cite{5} corresponding to $\mathfrak{R}(
s) =\sigma $ we have an analogous expression to the Baez-Duarte
formula for the $c_{k}$ appearing in the expansion of ${\zeta (
s) }^{-1}$ \cite{2,3}:
\begin{equation}
k^{\frac{\alpha -\sigma }{\beta }}\hat{d_{k}}=\frac{1}{\beta }\sum
\limits_{\rho }\Gamma ( \frac{\alpha -\rho }{\beta }) k^{\frac{\rho
-\sigma }{\beta }}+\frac{1}{\beta }\sum \limits_{n=1}^{\infty }\Gamma
( \frac{\alpha +2n}{\beta }) k^{-\frac{2n+\sigma }{\beta }}=:\psi
_{1}( k) 
\end{equation}

On the other hand we can express ${\hat{d}}_{k}$ and then the critical
function with a second formula:
\begin{gather}
{\hat{d}}_{k}=\frac{1}{\beta }\Gamma ( \frac{\alpha -1}{\beta })
k^{-\frac{\alpha -1}{\beta }}-\sum \limits_{p\,\mathrm{prime}}\ln
p\sum \limits_{q=1}^{\infty }\frac{1}{p^{\alpha q}}{\left( 1-\frac{1}{p^{\beta
q}}\right) }^{k}
\end{gather}
\begin{equation}
k^{\frac{\alpha -\sigma }{\beta }}\hat{d_{k}}=\frac{1}{\beta }\Gamma
( \frac{\alpha -1}{\beta }) k^{\frac{1-\sigma }{\beta }}-k^{\frac{\alpha
-\sigma }{\beta }}\sum \limits_{p\,\mathrm{prime}}\ln  p\sum \limits_{q=1}^{\infty
}\frac{1}{p^{\alpha q}}{\left( 1-\frac{1}{p^{\beta q}}\right) }^{k}=:\psi
_{2}( k) 
\end{equation}

In fact (see above) the Pochhammer expansion for $\frac{1}{s-1}$
is:
\[
\frac{1}{s-1}=\sum \limits_{k=0}^{\infty }s_{k}P_{k}( \frac{s-\alpha
}{\beta }+1) 
\]

\noindent where
\[
s_{k}=\int _{0}^{\infty } e ^{-\lambda ( \alpha -1) }(1-{\left.
e ^{-\lambda \beta }\right) }^{k}d\lambda 
\]

\noindent which for large {\itshape k} behaves as $\frac{1}{\beta
}\Gamma ( \frac{\alpha -1}{\beta }) k^{-\frac{\alpha -1}{\beta }}$.
Indeed with the substitution $ e ^{-\lambda \beta }=x$ we obtain:
\[
s_{k}=\frac{1}{\beta }\int _{0}^{1}{x^{\frac{\alpha -1}{\beta }-1}(
1-x) }^{k}dx=\frac{1}{\beta }\int _{0}^{1}{x^{\frac{\alpha -1}{\beta
}-1}( 1-x) }^{k+1-1}dx=\frac{1}{\beta }B( \frac{\alpha -1}{\beta
},k+1) 
\]

It is interesting to note that one can express the critical function
in terms of the zeros of the Zeta function (3.8) or in terms of
the primes (3.10). We will investigate numerically these two functions
for the case $\alpha =\frac{9}{2}$, $\beta =4$, $\sigma =\frac{1}{2}$.

\section{Numerical experiments}

As a test of the goodness of (3.2) we draw in Figure 3 the plots
of the function $\ln ( (1-2^{1-\sigma })\zeta ( \sigma ) ) $ and
of its polynomial representation in the interval $\sigma \in [-1,1[$.
Figure 3 shows a good match between them also in the ``critical
real interval'' $[0,1]$. 
\begin{figure}[h]
\begin{center}
\includegraphics{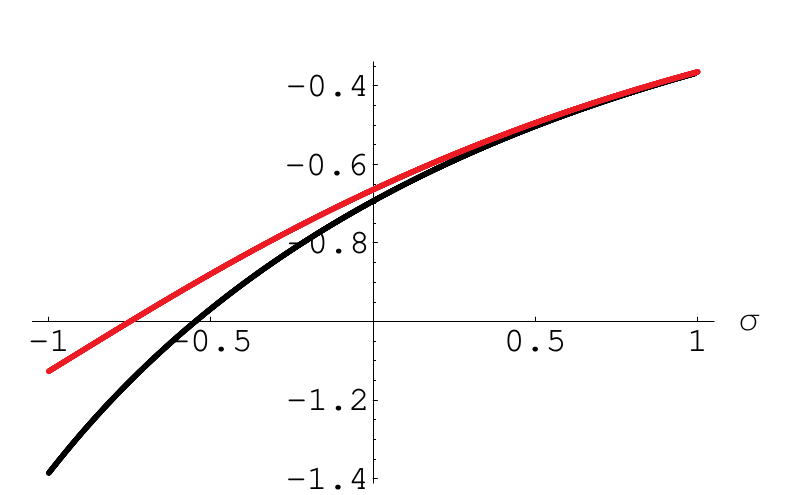}

\end{center}
\caption{The function $\ln ( (1-2^{1-\sigma })\zeta ( \sigma ) )
$ [black] and its polynomial representation [red]\ \ }

\end{figure}

In the next figures we present the results of the numerical experiment
performed on our representation (3.5) for the case $\alpha =\frac{9}{2}$
and $\beta =4$. We calculated the critical functions $\psi _{1}$
and $\psi _{2}$ for $\mathfrak{R}( z) =\sigma =\frac{1}{2}$. In
our calculations we considered only the first 10 non-trivial zeros
of the Zeta function, the first 20 trivial ones and the first 5'000
primes. Furthermore using the usual substitution $x=\log  k$, $\psi
_{1}$ and $\psi _{2}$ become:\ \ 
\begin{align*}
\psi _{1}( x) &=\frac{1}{4}\left( \sum \limits_{j=1}^{10}\Gamma
( 1-\frac{i t_{j}}{4})  e ^{\frac{xit_{j}}{4}}+\sum \limits_{j=1}^{10}\Gamma
( 1+\frac{i t_{j}}{4})  e ^{-\frac{xit_{j}}{4}}+\sum \limits_{n=1}^{20}\Gamma
( \frac{1}{2}n+\frac{9}{8})  e ^{-x( \frac{1}{2}n+\frac{1}{8}) }\right)
\\%
\psi _{2}( x) &=\frac{1}{4}\Gamma ( \frac{7}{8})  e ^{\frac{x}{8}}-
e ^{x}\sum \limits_{\overset{5000}{\mathrm{primes}}}\ln  p\sum \limits_{q=1}^{50}p^{-\frac{9}{2}q}
e ^{-\frac{ e ^{x}}{p^{4q}}}
\end{align*}

\noindent where $t_{j}$ is the imaginary part of the {\itshape j}-{\itshape
th} non-trivial zero.

We argue $\psi _{2}$ should approach $\psi _{1}$. The convergence
is surprising. The computations presented in Figure 4 and Figure
5 indicate that the qualitative and quantitative agreement between
the two functions is very good in the range $2.5\leq x\leq 30$ ($15\leq
k\leq 1.068\times {10}^{13}$).\ \ \ 
\begin{figure}[h]
\begin{center}
\includegraphics{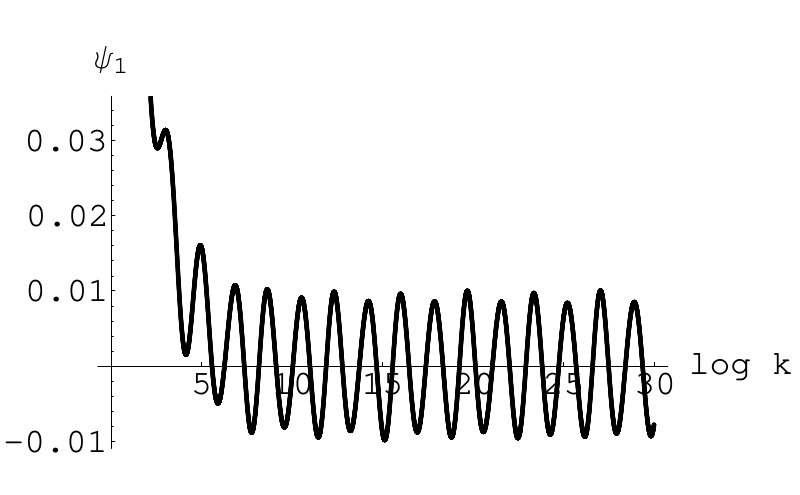}

\end{center}
\caption{The critical function $\psi _{1}$ calculated with the zeros
of the Zeta function}

\end{figure}
\begin{figure}[h]
\begin{center}
\includegraphics{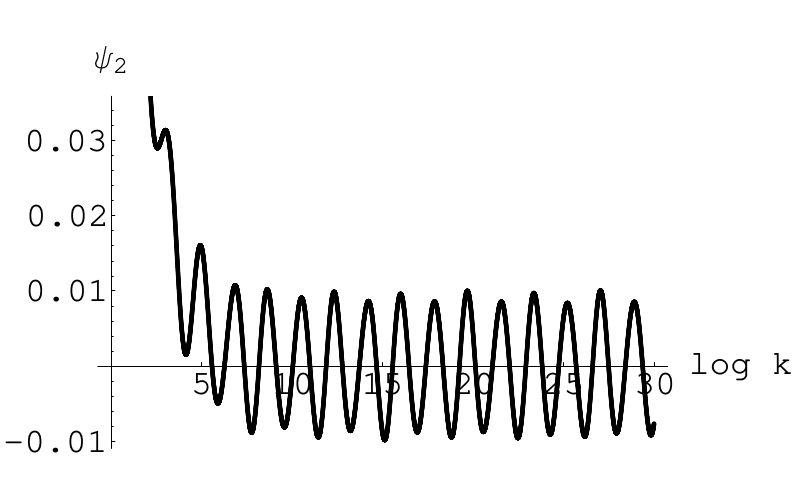}

\end{center}
\caption{The critical function $\psi _{2}$ calculated with the primes}

\end{figure}

It is interesting to study the single contribution of a prime to
the critical function $\psi _{2}$. In Figure 6 we computed the contributions
of the 10th prime ($p=29$), of the 50th prime ($p=229$) and of the
100th prime ($p=541$), all the calculations were performed until
$q=100$. The computations indicate that not only the contributions
decrease with increasing {\itshape p} but also that great primes
have an influence only on big values of {\itshape k}.\ \ \ \ \ 
\begin{figure}[h]
\begin{center}
\includegraphics{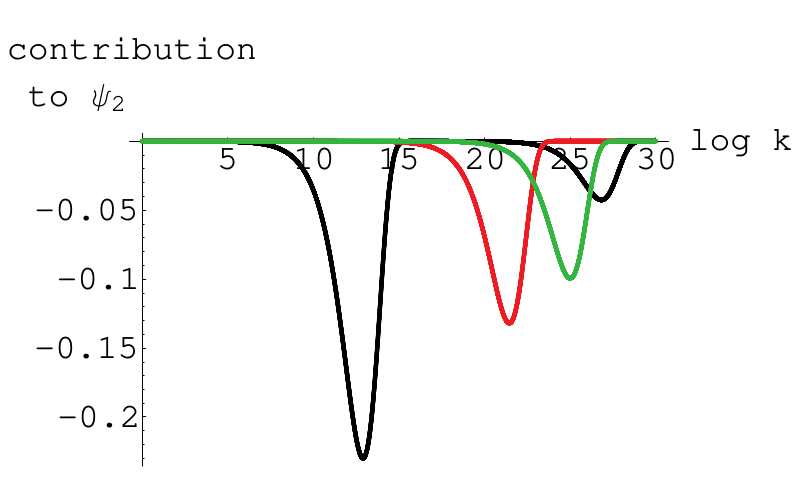}

\end{center}
\caption{The contribution to the critical function $\psi _{2}$ of
the primes $p=29$ [black], $p=229$ [red] and $p=541$ [green]}

\end{figure}

\section{Infinite $\beta$ limit}

In a numerical context we are also interested in the case of large
$\beta$ values. We start with the equation (7) in \cite{18}, given
by:
\begin{equation}
f( s) :=\left( \frac{\zeta ^{\prime }( s) }{\zeta ( s) }+\frac{1}{s-1}\right)
\frac{1}{s}=\sum \limits_{\rho }\frac{1}{\rho ( s-\rho ) }-\sum
\limits_{n=1}^{\infty }\frac{1}{2n( s+2n) }+\left( \frac{\zeta ^{\prime
}( 0) }{\zeta ( 0) }-1\right) \frac{1}{s}
\end{equation}

\noindent and we set $C=\frac{\zeta ^{\prime }( 0) }{\zeta ( 0)
}-1=\ln  2\pi -1$. Then, using the formula $\frac{1}{A}=\int _{o}^{\infty
} e ^{-\lambda  A}d\lambda $ as above ($\mathfrak{R}( A) >0 !$),
we obtain:
\begin{equation}
f( s) =\sum \limits_{k=0}^{\infty }{\hat{\hat{d}}}_{k}P_{k}( \frac{s-\alpha
}{\beta }+1) 
\end{equation}

\noindent where
\begin{equation}
{\hat{\hat{d}}}_{k}=\frac{1}{\beta }\left( \sum \limits_{\rho }\frac{\Gamma
( \frac{\alpha -\rho }{\beta }) }{\rho }k^{-\frac{\alpha -\rho }{\beta
}}-\sum \limits_{n=1}^{\infty }\frac{\Gamma ( \frac{\alpha +2n}{\beta
}) }{2n}k^{-\frac{\alpha +2n}{\beta }}+C \Gamma ( \frac{\alpha }{\beta
}) k^{-\frac{\alpha }{\beta }}\right) 
\end{equation}

We now analyze $\psi ( k) $, the absolute value of the critical
function, at large $\beta$ values where $\frac{1}{\beta }\Gamma
( \frac{\alpha -\rho }{\beta }) \sim \frac{1}{\alpha -\rho }$ is
valid.
\begin{equation}
\psi ( k) :=\frac{\left| {\hat{\hat{d}}}_{k}\right| }{k^{\frac{\sigma
-\alpha }{\beta }}}=\left| \sum \limits_{\rho }\frac{1}{\rho ( \alpha
-\rho ) }k^{-\frac{\sigma -\rho }{\beta }}-\sum \limits_{n=1}^{\infty
}\frac{1}{2n( \alpha +2n) }k^{-\frac{\sigma +2n}{\beta }}+\frac{C}{\alpha
}k^{-\frac{\sigma }{\beta }}\right| 
\end{equation}

Here the second and third term in the bracket converge for all $\sigma
>0$ (in particular for $\frac{1}{2}\leq \sigma \leq 1$). If we choose
$\alpha =1$, (5.4) would become in the $\beta$ limit (supposing
that this limit may be performed and has a meaning):
\begin{equation}
\operatorname*{\lim }\limits_{\beta \,\rightarrow \:\infty }\psi
( k) =\left| \sum \limits_{\rho }\frac{1}{\rho ( 1-\rho ) }-\sum
\limits_{n=1}^{\infty }\frac{1}{2n( 1+2n) }+C\right| \overset{\left(
5.1\right) }{=}\operatorname*{\lim }\limits_{x\,\rightarrow \:1}\left|
\frac{\frac{d}{dx}\ln  \zeta ( x)  +\frac{1}{x-1}}{x}\right| =\gamma
\end{equation}

\noindent where $\gamma \cong 0.577216$ is the Euler constant (see
also \cite{19}).

If such a limit is permitted our conjecture is that for $\mathfrak{R}(
s) \geq \sigma +\delta ,\delta >0$, as $\beta \overset{ }{\rightarrow
}\infty $: 
\begin{equation}
\left| f( s) \right| \sim \frac{B t}{\delta }\gamma 
\end{equation}

\noindent where {\itshape B} is some constant and $t=\mathfrak{I}(
s) $.

Since from the definition $P_{k}( \frac{s-\alpha }{\beta }+1) =\frac{\alpha
-s}{\beta }\frac{1}{k}P_{k-1}( \frac{s-\alpha }{\beta }) $ we obtain:
\[
\left| f( s) \right| \sim \sum \limits_{k=1}^{\infty }\left| {\hat{\hat{d}}}_{k}\right|
\left| \frac{1}{k}\frac{\alpha -s}{\beta }P_{k-1}( \frac{s-\alpha
}{\beta }) \right| 
\]

Then applying the Baez-Duarte inequality, i.e. $|P_{k-1}( z) |\leq
\frac{B}{{(k-1)}^{\mathfrak{R}( z) }}$ we have for $\mathfrak{R}(
s) \geq \sigma +\delta ,\delta >0$ that:
\[
\left| f( s) \right| \sim \sum \limits_{k=2}^{\infty }\frac{B\left|
\frac{\alpha -s}{\beta }\right| }{{k( k-1) }^{\frac{\sigma +\delta
-\alpha }{\beta }}}\left| {\hat{\hat{d}}}_{k}\right| \sim \sum \limits_{k=2}^{\infty
}\frac{B\left| \frac{\alpha -s}{\beta }\right| }{k^{\frac{\delta
}{\beta }+1}}\psi ( k) \sim B\frac{\left| \alpha -s\right| }{\beta
}\sum \limits_{k=2}^{\infty }\frac{1}{k^{1+\frac{\delta }{\beta
}}}\psi ( k) 
\]

\noindent and finally:
\[
\left| f( s) \right| \sim \left( B\frac{\left| \alpha -s\right|
}{\beta }\frac{{\left( k-1\right) }^{-\frac{\delta }{\beta }}}{-\frac{\delta
}{\beta }}\overset{\infty }{\operatorname*{|}\limits_{2}}\right)
\psi ( k) \sim B\frac{\left| \alpha -s\right| }{\beta }\frac{\beta
}{\delta }\psi ( k) \sim \frac{B\sqrt{{\left( \alpha -\sigma \right)
}^{2}+t^{2}}}{\delta }\psi ( k) 
\]

A similar (of course not rigorous) limit is formally obtained for
$\psi _{2}( k) $ using the primes along the lines for (3.4) to (3.10),
which, as $\beta$$ \rightarrow $$ \infty $ is given by:
\[
\operatorname*{\lim }\limits_{\beta \,\rightarrow \:\infty }k^{\frac{\alpha
-\sigma }{\beta }}\hat{d_{k}}=\frac{1}{\alpha }\frac{1}{\alpha -1}-\frac{1}{\alpha
}\sum \limits_{p \mathrm{prime}}\ln  p\sum \limits_{q=1}^{\infty
}\frac{1}{p^{\alpha q}}
\]

\noindent and thus \cite{19}:
\[
\operatorname*{\lim }\limits_{\alpha \,\rightarrow \:1^{+}}\frac{1}{\alpha
}\left( \frac{1}{\alpha -1}-\sum \limits_{p\,\mathrm{prime}}\ln
p\sum \limits_{q=1}^{\infty }\frac{1}{p^{\alpha q}}\right) =\operatorname*{\lim
}\limits_{\alpha \,\rightarrow \:1^{+}}\frac{1}{\alpha }\frac{d}{d\alpha
}\log ( \left( \alpha -1\right) \zeta ( \alpha ) ) =\gamma 
\]

We carried out some numerical experiments restricted to large $\beta$
values (until $\beta ={10}^{6}$), using the first 3600 known zeros
\cite{20}. The computations in Figure 7 indicate that for a fixed
{\itshape k}, within the limit of accuracy of our computations,
the difference between (5.4) and $\gamma$ approximately stabilizes
to less than 0.001 indipendently from the choice of {\itshape k}.
The difference is largely due only to the term involving the non-trivial
zeros. That is if we need a higher precision we have to consider
more non-trivial zeros in (5.4). 
\begin{figure}[h]
\begin{center}
\includegraphics{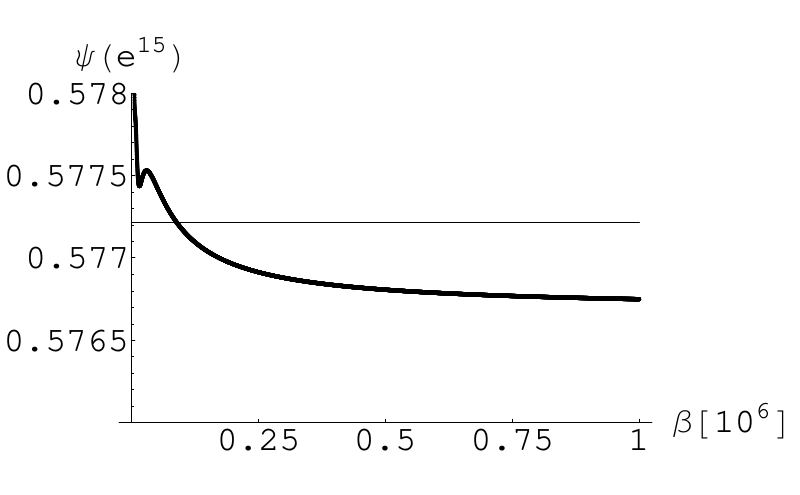}

\end{center}
\caption{The convergence to $\gamma$ of $\psi ( k) $ for $k= e ^{15}$}

\end{figure}

\section{Conclusions}

In this work we have found some new representations of functions
related to the Riemann Zeta function in terms of the Pochhammer's
polynomials, i.e. for the Zeta function via the alternating series,
for $(1-2^{1-s})\zeta ( s) $, for $\ln ( (1-2^{1-s})\zeta ( s) )
$ and for the derivative of $\ln ( (s-1)\zeta ( s) ) $.
\begin{enumerate}
\item A numerical experiment for the first function give satisfactory
results both for the real part as well for the imaginary part even
on the critical line $\mathfrak{R}( s) =\frac{1}{2}$ (we have used
the values $\alpha =\frac{1}{2}$, $\beta =i$ and {\itshape t} up
to $\mathfrak{I}( s) =t<35$).
\item In a formal limit of our representations (2.6) for the special
case $\alpha =\beta =2$ we obtain the Maslanka's representation
of $(s-1)\zeta ( s) $.
\item For the expansion of the derivative of the function $\ln (
(s-1)\zeta ( s) ) $ in terms of the Pochhammer's polynomials $P_{k}(
s) $ we have found two expressions ($\psi _{1}$ and $\psi _{2}$)
for the so called critical function: $\psi _{1}$ in terms of the
primes and $\psi _{2}$ in terms of the trivial as well as the non-trivial
zeros. We have then carried out a numerical experiment which gives
a very satisfactory agreements between the two, which up to very
high values of {\itshape k} remain bounded. The existence of absolut
upper bounds for the critical functions at {\itshape k}-infinity
may be considered as being equivalent to the truth of the RH.
\item Concerning the critical function in the large $\beta$ limit,
using $\alpha =1$, we may conjecture that $\frac{1}{s}$ time the
derivative of $\ln ( (s-1)\zeta ( s) ) $, using the inequality for
the Pochhammer's polynomials has, for $\mathfrak{R}( s) >\frac{1}{2}+\delta
$, a bound of the form $\frac{B t}{\delta }\gamma $ where $\gamma$
is the Euler constant and $t=\mathfrak{I}( s) $. 
\end{enumerate}

\end{document}